\documentclass{article}
\usepackage{amssymb,amsmath,amscd,amsthm}
\begin{document}

\newtheorem{theorem}{Theorem}
\newtheorem{lemma}{Lemma}
\newtheorem{definition}{Definition}
\newtheorem{remark}{Remark}

{\centerline {\Large {\bf Automorphisms  of the semigroup of invertible}}}

{\centerline {\Large {\bf  matrices with nonnegative elements}}}
\vspace{.3truecm}

\rightline{\large Bunina E.\,I., Mikhalev A.\,V.}

\vspace{.2truecm}

Let $R$ be a linearly ordered ring with $1/2$, $G_n(R)$ ($n\ge 3$) 
be a subsemigroup
of the group $GL_n(R)$, which consists of all
matrices with nonnegative elements. In the paper~[1]
A.V.~Mikhalev and M.A.~Shatalova described all automorphisms of
$G_n(R)$ in the case when  $R$ is a skewfield and $n\ge 2$. In this paper
we describe all automorphisms of $G_n(R)$, if
$R$ is an arbitrary linearly ordered associative ring
with $1/2$, $n\ge 3$.

\section{Main definitions and notations, formulation of the main theorem.}\leavevmode

Let $R$ be an associative ring with~$1$. 

\begin{definition}
A ring $R$ is called linearly ordered, if some subset $R_+$ 
is selected in~$R$, and this subset $R_+$ satisfies the following conditions\emph{:}

\emph{1)} $ \forall x\in R\ (x=0\lor x\in R_+\lor -x\in R_+)\land (x\in R_+
\Rightarrow -x\notin R_+)$\emph{;}

\emph{2)} $\forall x,y\in R_+ \ (x+y\in R_+\land xy\in R_+)$.

We  say that  $x\in R$ is greater \emph{(}not smaller\emph{)} 
than~$y\in R$, and  denote it by
 $x> y$, or $y< x$ \emph{(}$x\ge y$, or $y\le x$\emph{)}, iff $x-y\in R_+$
\emph{(}$x-y\in R_+\cup \{ 0\}$\emph{)}. 
\end{definition}

It is clear that $1\in R_+$, because
in the opposite case $-1\in R_+\Rightarrow 1=(-1)(-1)\in R_+$,
and that is impossible.

It is easy to prove by induction that in linearly ordered ring~$R$
we have $char\, R=0$.

Elements of  $R_+$ are called \emph{positive}, and elements
of $R_+\cup \{ 0\}$ are called \emph{nonnegative}.

\begin{definition}
Let $R$ be a linearly ordered ring. By  $G_n(R)$ we will denote
the subsemigroup of $GL_n(R)$, which consists of all matrices
with nonnegative elements.
\end{definition}

The set of all invertible elements of~$R$ is denoted by~$R^*$.
If $1/2\in R$, then $R^*$ is infinite, because it contains
all $1/2^n$ for $n\in \mathbb N$. The set $R_+\cap  R^*$ is denoted
by~$R_+^*$. If $1/2\in R$, then it is also infinite.

\begin{definition}
Suppose that $R$ is a linearly ordered ring, $T\subset R$. Then
 $Z(T)$ will denote  the center of~$T$, $Z^*(T)=Z(T)\cap R^*$,
$Z_+(T)=Z(T)\cap R_+$, $Z_+^*(T)=Z(T)\cap R_+^*$.
\end{definition}

It is clear that $Z_+^*(R)\subseteq Z_+^*(R^*)$. If $1/2\in R$, then
all these sets are infinite for $T=R$.

\begin{definition}
Let $I=I_n$, $\Gamma_n(R)$ be the group of all invertible
elements of $G_n(R)$, $\Sigma_n$ be the symmetric group of the order~$n$,
$S_\sigma$ be the matrix of a permutation $\sigma\in \Sigma_n$ \emph{(}i.\,e. 
the matrix $(\delta_{i\sigma(j)})$, where $\delta_{i\sigma(j)}$ is a Kronecker 
symbol\emph{)},
$S_n=\{ S_\sigma|\sigma\in \Sigma_n\}$, $diag[d_1,\dots,d_n]$ be
a diagonal matrix with the elements $d_1,\dots,d_n$ on the diagonal,
$d_1,\dots,d_n\in R_+^*$.
\end{definition}

\begin{definition}
By $D_n(R)$ we will denote the group of all invertible diagonal matrices
from $G_n(R)$, by $D_n^Z(R)$ we will denote the center of $D_n(R)$.
\end{definition}

It is clear that the group $D_n^Z(R)$ consists of all matrices $diag[d_1,\dots,
d_n]$, $d_1,\dots,d_n\in Z_+^*(R^*)$.

\begin{definition}
If  ${\cal A},{\cal B}$ are subsets of $G_n(R)$, then we will set
$$
C_{\cal A}({\cal B})=\{ a\in {\cal A}|\forall b\in {\cal B}\ (ab=ba)\}.
$$
\end{definition}

A matrix $A\in \Gamma_n(R)$, which satisfies $A^2=I$, will be called
an \emph{involution}.

\begin{definition}
By $K_n(R)$ we will denote  the subsemigroup of $G_n(R)$ consisting 
of all matrices
$$
\begin{pmatrix}
X_{n-1} & 0\\
0& x
\end{pmatrix},\quad X_{n-1}\in G_{n-1}(R),\ x\in R_+^*.
$$
\end{definition}

Let $E_{ij}$ be the matrix with one nonzero element $e_{ij}=1$.

\begin{definition}
By $B_{ij}(x)$ we will denote the matrix $I+xE_{ij}$.
Let ${\mathbf P}$ denote the subsemigroup in $G_n(R)$, which is
generated by the matrices $S_\sigma$ \emph{(}$\sigma\in \Sigma_n$\emph{)},
$B_{ij}(x)$ \emph{(}$x\in R_+, i\ne j$\emph{)}
 and $diag[\alpha_1,\dots,\alpha_n]\in D_n(R)$.
\end{definition}

\begin{definition}
Two matrices $A,B\in G_n(R)$ are called $\cal P$-equivalent
\emph{(}see~\emph{[1])},
if there exist matrices $A_j\in G_n(R)$, $j=0,\dots,k$, $A=A_0,B=A_k$, 
and matrices $P_i,\widetilde P_i, Q_i,\widetilde Q_i\in \mathbf P$, $i=0,\dots ,
k-1$ such that $P_iA_i\widetilde P_i=Q_iA_{i+1}\widetilde Q_i$.
\end{definition}

\begin{definition}
By $GE_n^+(R)$ we will denote the subsemigroup in $G_n(R)$ 
generated by all matrices
which are $\cal P$-equivalent to matrices from~$\mathbf P$.
\end{definition}

Note that if $R$ is a skewfield, then $GE_n^+(R)=G_n(R)$. 

\begin{definition}
If $G$ is some semigroup \emph{(}for example, $G=R_+^*, G_n(R), GE_n^+(R)$\emph{)}, 
then a 
homomorphism $\lambda(\cdot ): G\to G$ is called a central homomorphism
of~$G$, if $\lambda(G)\subset Z(G)$. The mapping $\Omega(\cdot): G\to G$ such that 
$\forall X\in G$
$$
\Omega (X)=\lambda(X)\cdot X,
$$
where $\lambda(\cdot)$ is a central homomorphism, is called a central
homothety.
\end{definition}

For example, if $R=\mathbb R$ (the field of real numbers), then the homomorphism
$\lambda(\cdot): G_n(\mathbb R)\to G_n(\mathbb R)$ such that $\forall A\in G_n(\mathbb 
R)$ $\lambda(A)=|det\, A|\cdot I$, is a central homomorphism, and the mapping
$\Omega(\cdot): G_n(\mathbb R)\to G_n(\mathbb R)$, such that $\forall A\in G_n(\mathbb R)$
$\Omega(A)=|det\, A|\cdot A$, is a central homothety.
Note that a central homothety $\Omega(\cdot)$ always is an endomorphism
of the semigroup~$G$: $\forall X,Y\in G$ $\Omega(X)\Omega(Y)=\lambda(X)X\cdot \lambda(Y)Y=
\lambda(X)\lambda(Y) X\cdot Y=\lambda(XY)XY=\Omega(XY).$

For every $M\in \Gamma_n(R)$ $\Phi_M$ denotes the automorphism 
of the semigroup $G_n(R)$ such that $\forall X\in G_n(R)$ 
$\Phi_M(X)=MXM^{-1}$.

For every $y(\cdot)\in Aut (R_+)$ by $\Phi^y$ we denote the automorphism 
of the semigroups $G_n(R)$ such that $\forall X=(x_{ij})\in G_n(R)$
$\Phi^y(X)=\Phi^y((x_{ij}))=(y(x_{ij}))$.

The main result of our paper is the following

{\bf Theorem.}
\emph{Let $\Phi$ be any automorphism of $G_n(R)$, $n\ge 3$, $1/2\in R$. 
Then on the semigroup $GE_n^+(R)$ 
$\Phi=\Phi_M\Phi^c\Omega$, where 
$M\in \Gamma_n(R)$, $c(\cdot)\in Aut(R_+)$, $\Omega(\cdot)$ is a central homothety of
$GE_n^+(R)$.}

\section{Construction of the automorphism $\Phi'$}\leavevmode

In this section we suppose that we have some automorphism $\Phi\in Aut(G_n(R))$,
where $n\ge 3$, $1/2\in R$, and we construct a new automorphism
$\Phi'\in Aut(G_n(R))$ such that $\Phi'=\Phi_{M'}\Phi$ for some $M'\in
\Gamma_n(R)$ and for all  $\sigma\in \Sigma_n$ we have $\Phi'(S_\sigma)=
S_\sigma$.

The following lemma in even a more general case can be found in~[2].

\begin{lemma} $\Gamma_n(R)=D_n(R)\cdot S_n$, i.e. the group $\Gamma_n(R)$ consists
of all monomial matrices.
\end{lemma}

\begin{proof}
It is clear that every monomial matrix is invertible, i.e.
$D_n(R)S_n\subset \Gamma_n(R)$.

Now let us consider some $A=(a_{ij})\in \Gamma_n(R)$.
 We need to show that in each its
row (column) there is exactly one nonzero element. Suppose that it is 
not true,
and the $i$-th row of~$A$ contains at least two nonzero (i.\,e. positive)
elements $a_{i,k}$ and $a_{i,j}$. Let us consider the inverse matrix 
$B=(b_{l,m})$.
Its  $k$-th row is nonzero, therefore there exists $l$ such that
 $b_{k,l}> 0$. Therefore
$$
\delta_{il}=a_{i,1}b_{1,l}+\dots+a_{i,n}b_{n,l}\ge a_{i,k}b_{k,l}>0,
$$
and so $i=l$. 

Similarly there exists $m$ such that $b_{j,m}> 0$, i.\,e.
$i=m$. Thus $l=m=i$.
Therefore $b_{j,i}> 0$, $b_{k,i}> 0$.

The condition $I=BA$ implies 
$$
\delta_{j,k} = b_{j,1}a_{1,k}+\dots+b_{j,n}a_{n,k} \ge b_{j,i} a_{i,k} > 0.
$$
Consequently $j=k$, which contradicts to the assumption that $i$-th row
contains two non-zero entries.
\end{proof}

Note that the representation of every matrix $A\in \Gamma_n(R)$ in the form
$A=D S_\sigma$, $D\in D_n(R)$, $\sigma\in \Sigma_n$, is unique.

\begin{lemma}
If $r\in R_+$ and $r^k=1$ for some $k\ge 1$, then $r=1$. 
\end{lemma}

\begin{proof}
 We need to show that  $x> 1$ implies $x^k> 1$,
and $0< x< 1$ implies $0< x^k< 1$.

1) Let us prove  by induction that
$$
x> 1\Rightarrow x^k> 1.
$$
If $k=1$, the condition is clear. Assume that the condition is proved for
some~$k$,
i.\,e.  $x> 1$, $x^k> 1$, therefore $x-1\in R_+$, $x^k-1\in R_+$, 
consequently, $x^{k+1}-x\in R_+\Rightarrow (x^{k+1}-x)+(x-1)\in R_+
\Rightarrow x^{k+1}-1\in R_+\Rightarrow x^{k+1}> 1$.

2) Similarly let us prove by induction that
$$
0< x< 1\Rightarrow 0< x^k< 1.
$$
If $k=1$, then the condition is clear. Assume that the condition is
 proved for some~$k$, i.\,e.
$0< x< 1\Rightarrow 0< x^k< 1$. Thus $x,x^k, 1-x, 1-x^k\in R_+$, 
therefore $x(1-x^k)=x-x^{k+1}\in R_+\Rightarrow (1-x)+(x-x^{k+1})=1-x^{k+1}
\in R_+$, i.\,e. $0< x^{k+1}< 1$.

The condition is proved.
\end{proof}

It is clear that in the ring $R$ there is no zero divizors.

The proof of the following lemma can be found in~[1].

\begin{lemma}
If $A$ is an involution in $G_n(R)$, then $A=diag [t_1,\dots,t_n]S_\sigma$, 
where $\sigma^2=1$, and for every $i=1,\dots, n$ 
we have $t_i\cdot t_{\sigma(i)}=1$.
\end{lemma}

\begin{proof} 
By Lemma~1 we have $A=dS_\sigma$, where $d=diag[d_1,\dots,d_n]$. Since
$A^2=I$, we have $dS_\sigma d S_\sigma=I \Rightarrow dS_\sigma =S_{\sigma}^{-1}
d^{-1} S_\sigma S_{\sigma}^{-1}$. Since the representation  of~$A$ in the form
$dS_\sigma$ is unique and $S_\sigma^{-1} d^{-1} S_\sigma\in D_n(R)$,
 we have  $d=S_\sigma^{-1}d^{-1}S_\sigma$ and $S_\sigma=S_\sigma^{-1}$.

Therefore $\sigma^2=1$ and $diag [d_1,\dots,d_n]=diag [d^{-1}_{\sigma(1)},\dots
, d^{-1}_{\sigma(n)}]$, i.\,e. $t_i=t_{\sigma(i)}^{-1}$.
\end{proof}

\begin{lemma}
If $\Phi$ is an automorphism of $G_n(R)$, where $n\ge 3$, $1/2\in R$,
then  

\emph{1)} $\Phi(\Gamma_n(R))=\Gamma_n(R),$

\emph{2)} $\Phi(D_n(R))=D_n(R),$

\emph{3)} $\Phi(D_n^Z(R))=D_n^Z(R).$
\end{lemma}

\begin{proof}
1) Since $\Gamma_n(R)$ is the subgroup of all invertible matrices of 
$G_n(R)$, then $\Phi(\Gamma_n(R))=\Gamma_n(R)$.

2) Consider the set $\cal F$ of all matrices  $A\in \Gamma_n(R)$, 
commuting with all matrices conjugate to~$A$.

Consider
$$
A=diag [\alpha_1,\dots,\alpha_n]\in D_n^Z (R),
$$
then every matrix conjugate to~$A$ has the form
\begin{multline*}
S_{\sigma^{-1}} diag[d_1^{-1},\dots,d_n^{-1}]diag [\alpha_1,\dots,\alpha_n]
diag [d_1,\dots,d_n] S_\sigma =\\
=S_{\sigma^{-1}}diag [\alpha_1,\dots,\alpha_n] S_\sigma =
diag [\alpha_{\sigma^{-1}(1)},\dots,\alpha_{\sigma^{-1}(n)}],
\end{multline*}
i.\,e. it commutes with $A=diag [\alpha_1,\dots,\alpha_n]$.

If we consider
$$
A=diag [\alpha_1,\dots,\alpha_n]\in D_n(R)\setminus D_n^Z(R),
$$
 then the matrix conjugate to~$A$ also is diagonal, but it is 
impossible to say, if
it commutes with
$A=diag [\alpha_1,\dots,\alpha_n]$, or not. 

Now consider
$$
A=diag [\alpha_1,\dots,\alpha_n]S_\rho,\quad \rho\ne e,\ \alpha_1,\dots,\alpha_n
\in R_+^*.
$$
 
Consider some 
$$
M=diag [d_1,\dots,d_n]\in D_n^Z(R),
$$
to obtain a matrix conjugate to~$A$.

We have 
\begin{multline*}
M^{-1}AM=diag[d_1^{-1},\dots,d_n^{-1}]diag [\alpha_1,\dots,\alpha_n] S_\rho
diag [d_1,\dots,d_n]=\\
=diag [d_{\rho^{-1}(1)}d_1^{-1}\alpha_1,\dots,d_{\rho^{-1}(n)}d_n^{-1}
\alpha_n] S_\rho=diag[\gamma_1\alpha_1,\dots,\gamma_n\alpha_n]S_\rho,
\end{multline*}
where $\gamma_1,\dots,\gamma_n\in Z_+^*(R^*)$.

We have the conditions 
\begin{align*}
A(M^{-1} AM)& = diag [\gamma_{\rho^{-1}(1)} \alpha_1 \alpha_{\rho^{-1}(1)},
\dots, \gamma_{\rho^{-1}(n)}\alpha_n \alpha_{\rho^{-1}(n)} ]S_{\rho^2},\\
(M^{-1} AM)A& = diag [\gamma_1 \alpha_1 \alpha_{\rho^{-1}(1)},\dots,
\gamma_n \alpha_n \alpha_{\rho^{-1}(1)}] S_{\rho^2}.
\end{align*}
Since $\rho\ne e$, we have that there exists such $i\in \{ 1,\dots,n\}$
that $j=\rho^{-1}(i)\ne i$. In this case we can take
$$
d_k=\begin{cases}
2,& \text{ if } k=j,\\
1,& \text{ if } k\ne j .
\end{cases}
$$
Thus we have 
$$
\gamma_k= \begin{cases}
2& \text{ if } k=i,\\
1/2,&\text{ if } k=j,\\
1,& \text{ if } k\ne i\text{ and } k\ne j.
\end{cases}
$$
Therefore
$$
A(M^{-1}AM)\ne (M^{-1}AM)A.
$$
So the condition
$$
(A\in \Gamma_n(R)) \land (\forall M\in \Gamma_n(R)\
(M^{-1}AM)A=A(M^{-1}AM))
$$
holds for all elements of $D_n^Z(R)$,  possibly 
can hold for some elements of $D_n(R)\setminus D_n^Z(R)$
and never holds for elements of $\Gamma_n(R)\setminus D_n(R)$.

It is clear that $\Phi({\cal F})={\cal F}$.

Introduce on the set $\cal F$ the additional condition
\begin{equation}
(A\in {\cal F})\land (\forall M\in \Gamma_n(R) \ 
(M\ne I\land M^{n!}=I\Rightarrow
AM\ne MA)),
\end{equation}
i.e. ``$A$ does not commute with any non-unit matrix of finite order''.

It is clear that if $A\in D_n(R)$ has two equal elements
on its $i$-th and $j$-th diagonal places, then it commutes with $S_{(i,j)}$.

If  $A\in D_n^Z(R)$ has all distinct entries  on its diagonal,
then it satisfies the condition~(1).

Besides, this condition can be satisfied for some 
 matrices from $D_n(R)\setminus D_n^Z(R)$, but they also must have
distinct entries on the diagonal. Let us denote the set of all matrices
satisfying the condition~(1), by~$\cal L$. It is clear that $\Phi({\cal L})=
{\cal L}$.

Let us consider the set
$$
{\cal X}=\bigcup_{M\in \cal L} C_{\Gamma_n(R)}(M),
$$
i.\,e. the set of all invertible matrices commuting with some 
matrix from~$\cal L$. 

Let us show that  ${\cal X}=D_n(R)$.

To prove ${\cal X}\subset D_n(R)$ note that every $M\in \cal L$ has distinct 
values on its diagonal, therefore if $AM=MA$, then $A\in D_n(R)$.

To prove that $D_n(R)\subset {\cal X}$ consider the matrix
$$
M=diag [2,2^2,\dots,2^n]\in D_n^Z(R).
$$
It is clear that $M\in {\cal L}$ and $C_{\Gamma_n(R)} (M)=D_n(R)$. 
Since $C_{\Gamma_n(R)}(M)\subset {\cal X}$, we have $D_n({\cal X})\subset 
{\cal X}$.

Consequently ${\cal X}=D_n(R)$.

It is clear that $\Phi({\cal X})={\cal X}$. Therefore
 $\Phi(D_n(R))=D_n(R)$.

3) Since $C_{\Gamma_n(R)}(D_n(R))=D_n^Z(R)$ and
$\Phi(D_n(R))=D_n(R)$, we have $\Phi(D_n^Z(R))=D_n^Z(R)$.
\end{proof}

\begin{lemma} If $\Phi$ is an automorphism of $G_n(R)$, $n\ge 3$,
$1/2\in R$, then there exists $M\in \Gamma_n(R)$ such that
$\Phi_M \Phi(K_n(R)) =K_n(R)$, where for all $X\in G_n(R)$
$$
\Phi_M(X)=MXM^{-1}.
$$
\end{lemma}

\begin{proof}
Let us consider a matrix 
$$
A=diag[\alpha,\dots,\alpha,\beta]\in D_n^Z(R),\quad \alpha\ne \beta.
$$
Suppose that
$$
B=\Phi(A)=diag[\gamma_1,\dots,\gamma_n]\in D_n^Z(R).
$$
It is clear that
$$
C_{\Gamma_n(R)} (A)/ D_n(R)\cong \Sigma_{n-1},
$$
hence
$$
\Phi(C_{\Gamma_n(R)}(A))/\Phi(D_n(R))=
 C_{\Gamma_n(R)}(B)/D_n(R)\cong \Sigma_{n-1}.
$$
Therefore 
$$
B=diag [\gamma,\dots,\gamma,\delta,\gamma,\dots,\gamma],\quad \gamma\ne \delta.
$$

Thus there  exists  $\sigma\in \Sigma_n$ such that
$$
\widetilde B=S_\sigma B S_{\sigma^{-1}}=[\gamma,\dots,\gamma,\delta].
$$
We have  $C_{G_n(R)}(A)\subseteq K$ and
$C_{G_n(R)}(\widetilde B)\subseteq K$. Besides, there exists such
a matrix $A$ (for example, $diag[1,\dots,1,2])$, that $C_{G_n(R)}(A)=K_n(R)$, 
i.\,e.
\begin{align*}
K_n(R)&=\!\!\! \bigcup_{A=diag[\alpha,\dots,\alpha,\beta]\in D_n^Z(R), 
\alpha\ne \beta} \!\!\!
C_{G_n(R)}(A),\\
K_n(R)&= \!\!\! \bigcup_{A=diag[\alpha,\dots,\alpha,\beta]\in D_n^Z(R), 
\alpha\ne \beta} \!\!\!
C_{G_n(R)}( S_\sigma \Phi(A) S_{\sigma}^{-1}),
\end{align*}
since $\Phi$ is an automorphism.

Consider $M=S_\sigma$, $\Phi'=\Phi_M \circ \Phi$. For every 
$A=diag [\alpha,\dots,\alpha,\beta]\in D_n^Z(R)$, $\alpha\ne \beta$, we have
$\Phi'(A)=diag [\gamma,\dots,\gamma,\delta]\in D_n^Z(R)$, $\gamma\ne \delta$, and
\begin{multline*}
\Phi'(K_n(R))=\Phi'\left( \bigcup_{A=diag 
[\alpha,\dots,\alpha,\beta]\in D_n^Z(R),\alpha
\ne \beta}\!\!\! C_{G_n(R)}(A)\right)=\\
=\bigcup_{A=diag[\alpha,\dots,\alpha,\beta]\in D_n^Z(R),\alpha\ne \beta}\!\!\!
\Phi'(C_{G_n(R)}(A))=\bigcup_{A=diag[\alpha,\dots,\alpha,\beta]
\in D_n^Z(R),\alpha
\ne \beta}\!\!\! C_{G_n(R)}(\Phi'(A))=K_n(R).
\end{multline*}
Therefore $\Phi_M\Phi(K_n(R))=K_n(R)$.
\end{proof}

\begin{lemma}
If $\Phi$ is an automorphism of $G_n(R)$, 
then by Lemma~\emph{4}, since
$\Phi(\Gamma_n(R))=\Gamma_n(R)$, for every $\sigma\in \Sigma_n$ we have
$$
\Phi(S_\sigma)=D_\sigma S_{\varphi(\sigma)},
$$
where $D_\sigma\in D_n(R)$. The obtained mapping $\varphi: \Sigma_n\to 
\Sigma_n$ is an endomorphism of the group~$\Sigma_n$.
\end{lemma}
\begin{proof}
For all $\sigma_1,\sigma_2\in \Sigma_n$
\begin{align*}
\Phi(S_{\sigma_1}\cdot S_{\sigma_2})&=\Phi(S_{\sigma_1\cdot \sigma_2})=
D_{\sigma_1\sigma_2}\cdot S_{\varphi(\sigma_1\sigma_2)},\\
\Phi(S_{\sigma_1})\cdot \Phi(S_{\sigma_2})&=D_{\sigma_1} S_{\varphi(\sigma_1)}
D_{\sigma_2} S_{\varphi(\sigma_2)}=D_{\sigma_1} \cdot D_{\sigma_2}'
S_{\varphi(\sigma_1)\varphi(\sigma_2)},\\
\Phi(S_{\sigma_1}\cdot S_{\sigma_2})&=\Phi 
(S_{\sigma_1})\cdot \Phi(S_{\sigma_2}),
\end{align*}
therefore 
$$
\varphi(\sigma_1\sigma_2)=\varphi(\sigma_1)\varphi(\sigma_2).
$$
Since $\Phi(D_n(R))=D_n(R)$ (see Lemma~4(2)), we have that 
$\varphi\in Aut(\Sigma_n)$. Actually, if $\sigma_1,\sigma_2\in \Sigma_n$, $\sigma_1
\ne \sigma_2$, $\sigma=\varphi(\sigma_1)=\varphi(\sigma_2)$, then
$$
\Phi(S_{\sigma_1})=D_{\sigma_1} S_\sigma,\
\Phi(S_{\sigma_2})=D_{\sigma_2} S_\sigma
\Rightarrow 
\Phi (S_{\sigma_1}\cdot S_{\sigma_2^{-1}})=D_{\sigma_1} S_\sigma S_{\sigma^{-1}}
D_{\sigma_2}^{-1}\in D_n(R).
$$
So for some $\rho\ne e$ we have
$$
\Phi(S_\rho)\in D_n(R),
$$
but it is impossible.
\end{proof}

The proof of the following lemma is completely similar to the proof of
Proposition~10 from~[1].

\begin{lemma}
If $\Phi$ is an automorphism of $G_n(R)$, $1/2\in R$, $n\ge 3$,
then there exists  a matrix
$M\in \Gamma_n(R)$ such that $\Phi_M \Phi(S_\sigma)=S_\sigma$
for all $\sigma\in \Sigma_n$.
\end{lemma}

\begin{proof}
1) Let $n\ne 6$. Consider the automorphism~$\varphi\in Aut(\Sigma_n)$ 
introduced in lemma~6. 
Since for $n\ne 6$ all automorphisms of the group $\Sigma_n$ are inner, 
then there exists a permutation $\rho \in \Sigma_n$ such that
for all $\sigma \in \Sigma_n$
$$
\varphi(\sigma)=\rho \sigma \rho^{-1}.
$$
Therefore for every $\sigma\in \Sigma_n$
$$
\Phi(S_\sigma)=D_\sigma S_{\rho \sigma \rho^{-1}} =D_{\sigma} S_\rho
S_\sigma S_{\rho^{-1}} =S_\rho D_\sigma' S_\sigma S_\rho.
$$
Then for $M'=S_\rho^{-1}$
$$
\forall \sigma \in \Sigma_n\ \ \Phi_{M'} \Phi (S_\sigma)=D_\sigma'S_\sigma.
$$

2) Let $n=6$. Consider some automorphism $\Phi_1=\Phi_{M_1}\Phi$ of
$G_6(R)$ such that $\Phi_1(K_n(R))=K_n(R)$, which exists by Lemma~5. 
Let $\varphi_1$ be
the automorphism of~$\Sigma_6$ defined by~$\Phi_1$. Note that
for every $\sigma \in \Sigma_6$ the matrix $S_\sigma\in  K_n(R)$ 
iff $\sigma(6)=6$.
Therefore $\varphi_1$ induces an inner automorphism of the group
$\Sigma (1,\dots,5)$, i.e. there exists $\tau \in \Sigma(1,\dots,5)$
such that for every $\sigma\in \Sigma(1,\dots,5)$ we have $\varphi_1(\sigma)=
\tau \sigma\tau^{-1}$.

Consider an automorphism
$\Phi_2=\Phi_{M_2}\circ \Phi_1$ of $G_6(R)$, where $M_2=S_\tau^{-1}$.
For every $\sigma\in \Sigma(1,\dots,5)$ we have
$$
\Phi_2(S_\sigma)=S_\tau^{-1} \Phi_1 (S_\sigma)S_\tau =S_\tau^{-1} D_\sigma
S_{\tau \sigma \tau^{-1}} S_\tau=D_\sigma' S_\sigma.
$$
Let $\varphi_2$ be an automorphism of~$\Sigma_6$ associated with~$\Phi_2$.
Note that $\Phi_2(K_n(R))=K_n(R)$ and $\varphi_2(\sigma)=\sigma$ 
for all $\sigma\in \Sigma (1,\dots,5)$. 

Let us prove that $\varphi_2$ is an identical
automorphism of~$\Sigma_6$. Let $\delta=\varphi_2((1,6))$. It is clear that
$\delta$ is an odd substitution, and either $\delta=(1,6)$, or
$\delta =(1,6)(\alpha_1,\alpha_2)(\beta_1,\beta_2)$. If $\delta=(1,6)$,
then for every $i=2,\dots,5$ 
$$
\varphi_2((i,6))=\varphi_2((1,6)(1,i)(1,6))=(1,6)(1,i)(1,6)=(i,6).
$$
Since the group $\Sigma_n$ is generated by transpositions and for every transposition~$
\sigma$ $\varphi_2(\sigma)=\sigma$, we have that for all $\sigma\in \Sigma_6$
$\varphi_2(\sigma)=\sigma$, i.e. $\varphi_2$ is identical.

 Let us show that the second
case is impossible. Actually, in this case $\varphi_2((1,6))=\varphi_2((\alpha_2,\beta_1)
(1,6)(\alpha_2,\beta_1))=(\alpha_2,\beta_1)(1,6)(\alpha_1,\alpha_2)(\beta_1,
\beta_2)(\alpha_2,\beta_1)=(1,6)(\alpha_1,\beta_1)(\alpha_2,\beta_2)\ne
\varphi(1,6)$.

Therefore for the matrix $M'=M_2M_1$ we have $\Phi_{M'} \Phi(S_\sigma)=D_\sigma 
S_\sigma$ for every $S_\sigma\in S_n$.

3) 
Now for every $n\ge 3$ we have some  $M'\in \Gamma_n(R)$ such
that $\Phi_{M'} \Phi(S_\sigma)=D_\sigma S_\sigma$ for every matrix
$S_\sigma \in S_n$. 

Consider  $\rho=(1,2,\dots,n)\in \Sigma_n$.
Let $\Phi_{M'} \Phi(S_\rho)=D_\rho S_\rho$, 
where $D_\rho =diag[\alpha_1,\dots,
\alpha_n]$. The equality  $S_\rho^n=I$ implies $(D_\rho S_\rho)^n=I$, 
and then $\alpha_1 \alpha_2\dots \alpha_n=1$. Consider a matrix
$T=diag[t_1,\dots,t_n]$, where $t_i=(\alpha_i\alpha_{i+1}\dots \alpha_n)^{-1}$,
$i=1,\dots,n$, and the automorphism $\Phi_3=\Phi_T\Phi_M\Phi$. We see that
$\Phi_3(S_\rho)=diag [t_1,\dots,t_n]diag [\alpha_1,\dots,\alpha_n] S_\rho
diag [t_1^{-1},\dots ,t_n^{-1}]=diag [t_1 \alpha_1t_2^{-1},t_2\alpha_2 t_3^{-1},\dots
.t_n\alpha_n t_1^{-1}]S_\rho=S_\rho$ and for every permutation
 $\sigma \in \Sigma_n$ 
we have $\Phi_3(S_\sigma)=\widetilde D_\sigma S_\sigma$, where
$\widetilde D_\sigma \in D_n(R)$. Let now $\tau =(1,2)$. 
Then $S_\tau$ is an involution. According to Lemma~3, we obtain
$\Phi_3(S_\tau)=\widetilde D_\tau S_\tau$, where 
$$
\widetilde D_\tau=
diag[\beta,\beta^{-1}, 1,\dots,1]\quad (\beta\in R_+^*).
$$
 The condition
$\rho=(n,n-1)(n-1,n-2)\dots(2,1)$ implies $S_\rho=S_\tau^{\rho^{n-2}}\dots S_\tau$,
where $S_\tau^\delta=S_\delta^{-1}S_\tau S_\delta$. Consequently,
$S_\rho =\Phi_3(S_\tau^{\rho^{n-2}}\dots S_\tau)=(D_\tau S_\tau)^{\rho^{n-2}}
\dots D_\tau S_\tau$. 

Comparing nonzero elements, we obtain
$\beta^n=1$, hence $\beta=1$ (see Lemma~2). Therefore $\Phi_3(S_\tau)=
S_\tau$. Since $\rho$ and $\tau$ generate $\Sigma_n$, we have  $\Phi_3(S_\sigma)
=S_\sigma$ for every $\sigma\in \Sigma_n$.
\end{proof}

\section{The action of $\Phi'$ on diagonal matrices.}\leavevmode

In the previous section
 by our automorphism $\Phi$ we have constructed  a new automorphism
$\Phi'=\Phi_M\Phi$ such that $\Phi' (S_\sigma)=S_\sigma$ 
for all $\sigma\in \Sigma_n$.
We will suppose that this automorphism $\Phi'$ is fixed.

\begin{lemma}
If $n\ge 3$, $1/2\in R$, the automorphism $\Phi' \in Aut(G_n(R))$ is
such that $\forall \sigma\in \Sigma_n$
$\Phi'(S_\sigma)=S_\sigma$, then for all $\alpha,\beta\in R_+^*$ we have
$$
\Phi' (diag[\alpha,\beta,\dots,\beta])=diag[\gamma,\delta,\dots,\delta],\quad
\gamma,\delta\in R_+^*.
$$
If $\alpha\ne \beta$, then $\gamma\ne \delta$. If $\alpha,\beta\in Z_+^{*}(R^*)$, 
then $\gamma,\delta\in Z_+^{*}(R^*)$.
\end{lemma}

\begin{proof}
By Lemma~4 
$$
\Phi'(diag [\alpha,\beta,\dots,\beta])=diag [\gamma_1,\dots,\gamma_n].
$$
Since $\Phi'(S_{(i,i+1)})=S_{(i,i+1)}$ for all $i=2,\dots,n-1$, then for all $i=
2,\dots,n-1$ we have
\begin{multline*}
\Phi'(diag [\alpha,\beta,\dots,\beta])\Phi' (S_{(i,i+1)})=\Phi'(S_{(i,i+1)})\Phi'(
diag[\alpha,\beta,\dots,\beta])\Rightarrow\\
diag [\gamma_1,\dots,\gamma_n]S_{(i,i+1)}=S_{(i,i+1)}diag [\gamma_1,\dots,\gamma_n]
\Rightarrow \gamma_i=\gamma_{i+1}.
\end{multline*}
Therefore $\gamma_2=\gamma_3=\dots=\gamma_{n-1}=\gamma_n$ and we can assume that
$$
\Phi'(diag[\alpha,\beta,\dots,\beta])=diag [\gamma,\delta,\dots,\delta].
$$
If $\alpha\ne \beta$, then
$$
diag [\alpha,\beta,\dots,\beta] S_{(1,2)}\ne S_{(1,2)} diag [\alpha,\beta,\dots,\beta]
\Rightarrow diag [\gamma,\delta,\dots,\delta]S_{(1,2)} \ne S_{(1,2)}
diag [\gamma,\delta,\dots,\delta]\Rightarrow \gamma\ne \delta.
$$
If $\alpha,\beta\in Z^*(R^*)$, then $diag[\alpha,\beta,\dots,\beta]\in D_n^Z(R)$,
and by Lemma 4(3) $diag[\gamma,\delta,\dots,\delta]\in D_n^Z(R^*)$,
thus $\gamma,\delta\in Z^*(R^*)$.
\end{proof}

\begin{lemma}
If $n\ge 3$, $1/2\in R$, the automorphism $\Phi'\in Aut(G_n(R)$ is such that $\forall
\sigma\in \Sigma_n$ $\Phi'(S_\sigma)=S_\sigma$, then
for all $X\in G_2(R)$ we have
$$
\Phi'\begin{pmatrix}
X& 0& \dots& 0\\
0& 1 & \dots& 0\\
\hdotsfor{2}& \ddots& \dots\\
0& \hdotsfor{2}& 1
\end{pmatrix}
=\begin{pmatrix}
Y& 0& \dots& 0\\
0& a & \dots& 0\\
\hdotsfor{2}& \ddots& \dots\\
0& \hdotsfor{2}& a
\end{pmatrix},\quad \text{ where } Y\in G_2(R), a\in Z_+^{*}(R^*).
$$
\end{lemma}
\begin{proof}
Similarly to the proof of Lemma~8 we can prove that for every 
$$
A=diag [\alpha,\alpha,\beta,\dots,\beta]\in D_n(R),\quad \alpha\ne \beta,
$$
we have
$$
\Phi'(A)=diag [\gamma,\gamma,\delta,\dots,\delta]\in D_n(R),\quad \gamma\ne \delta.
$$
Consider now the set $\cal L$ of all involutions of the form
$$
diag[\xi,\xi^{-1},1,\dots,1]S_{(1,2)},\quad \xi\in R_+^*.
$$
For every such an involution~$M$ we have 
$$
N=\Phi'(M)=\Phi'(diag [\xi,\xi^{-1},1,\dots,1])S_{(1,2)}
$$
and $N$ is an involution. By Lemma~3  
$$
N=diag[\eta,\eta^{-1},1,\dots,1]S_{(1,2)}.
$$
If $\xi\in Z^*_+(R^*)$, then $\eta\in Z_+^*(R^*)$, if $\xi \notin Z_+^*(R^*)$, then $\eta\notin 
Z_+^*(R^*)$. We see that $\Phi'({\cal L})=\cal L$. 

The set of all matrices of the form
$$
diag[\alpha,\alpha,\beta,\dots,\beta], \quad \alpha,\beta\in R_+^*,\ 
\alpha\ne \beta
$$
will be denoted by $\cal M$. We know that $\Phi'({\cal M})=\cal M$.
Therefore 
$$
\Phi'(C_{\cal M}{\cal L})=C_{\cal M}{\cal L},
$$
i.e. for every $\mu \in Z_+^*(R^*)$, $\eta\in R_+^*$ we have 
$$
\Phi'(diag [\mu,\mu,\eta,\dots,\eta])=diag [\mu',\mu',\eta',\dots,\eta'],
$$
where $\mu'\in Z_+^*(R^*)$, $\eta'\in R_+^*$ and if $\eta\in R_+^*\setminus
Z_+^*(R^*)$, then $\eta'\in R_+^*\setminus
Z_+^*(R^*)$. 

By $\cal Z$ we will denote the set of all matrices
$$
\alpha I=diag [\alpha,\dots,\alpha],\quad \alpha \in Z_+^*(R^*).
$$

It is clear that $\Phi'({\cal Z})=\cal Z$.

Consider some matrix
$$
A=\begin{pmatrix}
X& 0& & \\
0& a& & \\
& & \ddots& \\
& & & a
\end{pmatrix},\quad  X\in G_2(R), a\in Z_+^{*}(R^*).
$$
This matrix satisfies the condition
\begin{equation}                                                                                   
\forall M\in C_{\cal M} {\cal L}\ \exists N\in {\cal Z}\ A(MN)=(MN)A
\land AS_{(3,4)}=S_{(3,4)}A\land\dots \land AS_{(n-1,n)} =S_{(n-1,n)}A.
\end{equation}
Actually, every $M\in C_{\cal M}{\cal L}$ has the form 
$$
M=diag [\mu,\mu,\eta,\dots,\eta],\quad \mu\in Z_+^*(R^*), \eta\in R^*_+.
$$
If $MA=AM$, we can take $N=I$. If $MA\ne AM$, i.e. $\mu\in Z_+^*(R^*)
\setminus Z_+^*(R)$ and 
$Xdiag [\mu,\mu]\ne diag [\mu,\mu]X$, then we can take 
$$
N=diag [\mu^{-1},\dots,\mu^{-1}]\in {\cal Z}.
$$
Then $MN=diag[1,1,\eta\mu^{-1},\dots,\eta\mu^{-1}]$ and $A(MN)=(MN)A$.

If some matrix $A$ satisfies the condition~(2), then the part 
$$
AS_{(3,4)}=S_{(3,4)}A\land\dots\land AS_{(n-1,n)}=S_{(n-1,n)}A
$$
implies
$$
A=\begin{pmatrix}
X& 0& & \\
0& a& & \\
& & \ddots& \\
& & & a
\end{pmatrix},\quad  X\in G_2(R), a\in R^*_+.
$$
If $a\in R^*_+\setminus Z_+^*(R^*)$, then there exists $b\in R_+^*$ such that $ab\ne ba$,
therefore for
$$
M=diag [1,1,b,\dots,b]\in C_{\cal M}{\cal L}
$$
we have $MA\ne AM$. For every $N=diag[\alpha,\dots,\alpha]\in \cal Z$
we have $A(MN)\ne (MN)A$, since $ab\alpha\ne b\alpha a$.
Thus the matrix 
$$
A=\begin{pmatrix}
X& 0& & \\
0& a& & \\
& & \ddots& \\
& & & a
\end{pmatrix}
$$
with $a\in R_+^*\setminus Z_+^*(R^*)$ can not satisfy the condition~(2). 
So we have $a\in Z_+^*(R^*)$. 

Therefore we see that a matrix $A$ has the form 
$$
A=\begin{pmatrix}
X& 0& & \\
0& a& & \\
& & \ddots& \\
& & & a
\end{pmatrix},\quad  X\in G_2(R), a\in Z_+^*(R^*)
$$
iff it satisfies the condition~(2). 

Since $\Phi'(S_{(i,i+1)})=S_{(i,i+1)}$ for all $i=3,\dots,n-1$, $\Phi'({\cal Z})
=\cal Z$, $\Phi'(C_{\cal M}{\cal L})=C_{\cal M}{\cal M}$, we 
obtain that if a matrix $A$ satisfies~(2), then the matrix $\Phi'(A)$
satisfies~(2).

Consequently, for $X\in G_2(R)$, $a\in Z_+^*(R^*)$ we have 
$$
\Phi'\begin{pmatrix}
X& 0& & \\
0& a& & \\
& & \ddots& \\
& & & a
\end{pmatrix}=
\begin{pmatrix}
Y& 0& & \\
0& b& & \\
& & \ddots& \\
& & & b
\end{pmatrix},\quad  Y\in G_2(R), b\in Z_+^*(R^*).
$$
 
\end{proof}

\begin{lemma} If $n\ge 3$, $1/2\in R$, the automorphism 
$\Phi'\in Aut(G_n(R))$ is such that
$\forall \sigma\in \Sigma_n$ $\Phi'(S_\sigma)=S_\sigma$, then
 for every $x\in Z_+^*
(R)$
$$
\Phi'(diag[x,1,\dots,1])=diag [\xi,\eta,\dots,\eta],\quad \xi,\eta\in Z_+^*(R).
$$
\end{lemma}
\begin{proof}
Since $x\in Z_+^*(R)$, we have that $x\in Z_+^*(R^*)$, therefore $A=diag
[x,1,\dots,1]\in D_n^Z(R)$, then by Lemma~8 $A'=\Phi'(A)=diag[\xi,\eta,\dots,\eta]$,
where $\xi,\eta\in Z_+^*(R^*)$.

Let $\cal Y$ denote the set of all matrices of the form
$$
\begin{pmatrix}
 a& 0&\dots & 0\\
0  & \ddots& \hdotsfor{2}\\
\hdotsfor{2}  & a& 0\\
\hdotsfor{2} & 0& X
\end{pmatrix},\quad  X\in G_2(R), a\in Z_+^*(R^*).
$$
It is clear that 
$\Phi'({\cal Y})=\cal Y$ (the proof is completely similar to the proof of Lemma~9). 

Let $\overline{\cal Z} $ denote the center of the semigroup $G_n(R)$. It is clear that
$$
\overline{\cal Z} =\{ \alpha I|\alpha\in Z_+^*(R)\}.
$$
We have $\Phi'(\overline{\cal Z})=\overline {\cal Z}$.

Every matrix $A=diag[x,1,\dots,1]$, $x\in Z_+^*(R)$
satisfies the condition
\begin{equation}
\forall M\in {\cal Y}\ MA=AM.
\end{equation}
The matrix $A'=\Phi'(A)$ also satisfies the condition~(3), therefore
$$
\forall M\in {\cal Y} \quad M\cdot diag[\xi,\eta,\dots,\eta]=
diag[\xi,\eta,\dots,\eta]\cdot M,
$$
or
$$
\forall X\in G_2(R)\quad X\circ diag [\eta,\eta]=diag [\eta,\eta]\circ X,
$$
thus $\eta\in Z_+^*(R)$.

Now we need to prove that $\xi\in Z_+^*(R)$.

We have 
$$
\Phi'(diag[1,x,1,\dots,1])=\Phi'(S_{(1,2)}diag [x,1,\dots,1]S_{(1,2)})
=S_{(1,2)}diag [\xi,\eta,\dots,\eta]S_{(1,2)}=diag[\eta,\xi,\eta,\dots,\eta]
$$
and, similarly,
$$
\Phi'(diag [1,1,x,1,\dots,1])=diag [\eta,\eta,\xi,\eta,\dots,\eta],
\dots, \Phi'(diag [1,\dots,1,x])=diag [\eta,\dots,\eta,\xi].
$$

Consequently,
$$
\Phi'(x\cdot I)=\Phi'(diag[x,1,\dots,1]\cdot diag [1,x,\dots,1]\dots
 diag[1,\dots,1,x])=\xi\eta^{n-1}\cdot I.
$$
Since $x\in Z_+^*(R)$, we have that $\xi \eta^{n-1}\in Z_+^*(R)$. 
Since (as we just have proved) $\eta\in Z_+^*(R)$, we have that $\eta^{n-1}\in Z_+^*(R)$,
therefore $\xi\in Z_+^*(R)$, as we needed to prove.
\end{proof}

\begin{lemma}
If $n\ge 3$, $1/2\in R$, the automorphism 
$\Phi'\in Aut(G_n(R))$ is such that $\forall \sigma\in 
\Sigma_n$ $\Phi'(S_\sigma)=S_\sigma$, then for every $x_1,x_2\in Z_+^*(R)$ such that
$x_1\ne x_2$, 
\begin{align*}
\Phi'(A_1)&= \Phi'(diag [x_1,1,\dots,1])=diag [\xi_1,\eta_1,\dots,\eta_1],\\
\Phi'(A_2)&= \Phi'(diag [x_2,1,\dots,1])=
diag [\xi_2,\eta_2,\dots,\eta_2]
\end{align*}
we have $\xi_1\eta_1^{-1}\ne \xi_2\eta_2^{-1}$.
\end{lemma}
\begin{proof}
Suppose that for some distinct $x_1,x_2\in  Z_+^*(R)$ we have $\xi_1\eta_1^{-1}=
\xi_2\eta_2^{-1}$, i.e. 
\begin{align*}
\Phi'(A_1)&= \Phi'(diag [x_1,1,\dots,1])=diag [\xi,\eta,\dots,\eta]=A_1',\\
\Phi'(A_2)&= \Phi'(diag [x_2,1,\dots,1])=\alpha\cdot
 diag [\xi,\eta,\dots,\eta]=A_2',
\end{align*}
where $\xi,\eta,\alpha\in Z_+^*(R)$ (see Lemma~10). Therefore ${\Phi'}^{-1}
(\alpha I)={\Phi'}^{-1}(A_1' {A_2'}^{-1})=
diag[x_1x_2^{-1},1,\dots,1])=diag[\beta,1,\dots,1]$, where $1\ne \beta\in Z_+^*
(R)$, but it is impossible, since ${\Phi'}^{-1}(\overline{\cal Z})=\overline{\cal Z}$
(see the proof of Lemma~10).
Thus $\xi_1\eta_1^{-1}\ne \xi_2\eta_2^{-1}$.
\end{proof}

\section{The main theorem.}\leavevmode

In this section we will prove the main theorem (Theorem~1).

Recall (see Definition~8) that for $x\in R_+$
$$
B_{12}(x)=\begin{pmatrix}
1& x& &&\\
0& 1& &&\\
&& 1&&\\
&&&\ddots&\\
&&&& 1
\end{pmatrix},\quad
B_{21}(x)=\begin{pmatrix}
1& 0& &&\\
x& 1& &&\\
&& 1&&\\
&&&\ddots&\\
&&&& 1
\end{pmatrix}.
$$

\begin{lemma}
If $n\ge 3$, $1/2\in R$, the automorphism $\Phi'\in Aut(G_n(R))$ is such that 
$\forall \sigma\in \Sigma_n$ $\Phi'(S_\sigma)=S_\sigma$, then there are
two possibilities\emph{:}

\emph{1)} there exists some mapping $c(\cdot ): R_+\to R_+$ such that
for all $x\in R_+$ $\Phi'(B_{12}(x))=B_{12}(c(x))$\emph{;}

\emph{2)} there exists some mapping $b(\cdot ): R_+\to R_+$ such that
for all $x\in R_+$ $\Phi'(B_{12}(x))=B_{21}(b(x))$.
\end{lemma}

\begin{proof}
By Lemma~9 we have
$$
\Phi'(B_{12}(1))=\begin{pmatrix}
\alpha& \beta&&&\\
\gamma& \delta&&&\\
&&a&&\\
&&&\ddots&\\
&&&&a
\end{pmatrix}, \quad a\in Z_+^{*}(R^*), 
\begin{pmatrix}
\alpha& \beta\\
\gamma&\delta
\end{pmatrix} \in G_2(R)
$$
(see Lemma~9).

Let for every $x\in R_+^*$ 
$$
\Phi'(diag[x,1,\dots,1)]=diag[\xi(x),\gamma(x),\dots,\gamma(x)],\quad
\xi(x),\eta(x)\in R_+^*
$$
(see Lemma~8).

Then for every $x\in Z_+^*(R)$
\begin{multline*}
\Phi' (B_{12}(x))=\Phi' (diag[x,1,\dots,1]
B_{12}(1) diag[x^{-1},1,\dots,1])=\\
=diag [\xi(x),\eta(x),\dots,\eta(x)]
\begin{pmatrix}
\alpha& \beta& &&\\
\gamma& \delta&&&\\
&&a&&\\
&&&\ddots&\\
&&&&a
\end{pmatrix}
diag[\xi(x)^{-1},\eta(x)^{-1},\dots ,\eta(x)^{-1}]=\\
=\begin{pmatrix}
\xi(x)\alpha\xi(x)^{-1}& \xi(x)\beta\eta(x)^{-1}& &&\\
\eta(x)\gamma\xi(x)^{-1}& \eta(x)\delta\eta(x)^{-1}&&&\\
&&a&&\\
&&&\ddots&\\
&&&&a
\end{pmatrix}.
\end{multline*}
Since by Lemma~10 $\xi(x),\eta(x)\in Z_+^*(R)$, we have 
$\xi(x)\alpha \xi(x)^{-1}=\alpha$, $\xi(x)\beta \eta(x)^{-1}=\xi(x)\eta(x)^{-1}
\beta$, $\eta(x)\gamma \xi(x)^{-1}=\eta(x)\xi(x)^{-1}\gamma$, $\eta(x)\delta \eta(x)^{-1}=
\delta$, thus
$$
\Phi'(B_{12}(x))=
\begin{pmatrix}
\alpha& \nu(x)\beta& &&\\
\nu(x)^{-1} \gamma& \delta&&&\\
&&a&&\\
&&&\ddots&\\
&&&&a
\end{pmatrix}
$$
for $\nu(x)=\xi(x)\eta(x)^{-1}$.

By Lemma~11 for $x_1\ne x_2$ we have $\nu(x_1)\ne \nu(x_2)$. 

For every $x\in R_+$ $\Phi'(B_{12}(1))$ and $\Phi' (B_{12}(x))$ commute.
Let us write this condition in the matrix form for $x\in Z_+^*(R)$
(recall that in this case $\nu(x)\in Z_+^*(R)$ by Lemma~10):
\begin{multline*}
\begin{pmatrix}
\alpha& \beta\\
\gamma& \delta
\end{pmatrix}
\begin{pmatrix}
\alpha& \nu(x)\beta\\
\nu(x)^{-1}\gamma & \delta
\end{pmatrix} =\begin{pmatrix} \alpha& \nu(x)\beta\\
\nu(x)^{-1} \gamma& \delta
\end{pmatrix}
\begin{pmatrix}
\alpha& \beta\\
\gamma& \delta
\end{pmatrix} \Rightarrow\\
\Rightarrow
\begin{pmatrix}
\alpha^2+\nu(x)^{-1} \beta\gamma& \nu(x)\alpha\beta+\beta\delta\\
\gamma\alpha+\nu(x)^{-1}\delta\gamma& \nu(x)\gamma \beta+\delta^2
\end{pmatrix}=
\begin{pmatrix}
\alpha^2+\nu(x) \beta\gamma& \alpha\beta+\nu(x)\beta\delta\\
\nu(x)^{-1}\gamma\alpha+\delta\gamma& \nu(x)^{-1}\gamma \beta+\delta^2
\end{pmatrix}.
\end{multline*}
Hence $\nu(x)^{-1}\beta\gamma =\nu(x)\beta\gamma$ for distinct $x\in Z_+^*(R)$ 
(for example, for $x=2,2^2,\dots$). By Lemma~11 $\nu(x)\ne 1$ for 
$x\ne 1$, therefore $\nu(x)\ne \nu(x)^{-1}$ for $x\ne 1$ and $\beta\gamma=0$,
i.e. either $\beta=0$, or $\gamma=0$.

Suppose that $\gamma =0$ (the case $\beta=0$ is similar).

Then 
$$
\Phi'(B_{12}(1))=
\begin{pmatrix}
\alpha& \beta& &&\\
0& \delta&&&\\
&&a&&\\
&&&\ddots&\\
&&&&a
\end{pmatrix}, \quad a\in Z_+^*(R^*),\ \alpha,\delta\in R_+^*, \beta\in R_+\cup \{ 0\}.
$$
Let us use the condition $(B_{12}(1))^2=diag[2,1,\dots,1]B_{12}(1)\cdot diag [1/2,
1,\dots,1]$:
$$
\begin{pmatrix}
\alpha^2& \alpha\beta+\beta \delta& &&\\
0& \delta^2&&&\\
&&a^2&&\\
&&&\ddots&\\
&&&&a^2
\end{pmatrix}=
\begin{pmatrix}
\alpha& \nu(2)\beta& &&\\
0& \delta&&&\\
&&a&&\\
&&&\ddots&\\
&&&&a
\end{pmatrix},
$$
which implies $\alpha=\delta=a=1$, $\nu(2)=2$. Therefore we have $\Phi'(B_{12}
(1))=B_{12}(\beta)$ for some $\beta\in R_+$.

Similarly, if $\beta=0$, then $\Phi'(B_{12}(1))=B_{21}(\gamma)$ for some $\gamma\in 
R_+$.

Consider the case $\gamma=0$ (the case $\beta=0$ is similar). Since for every
$x\in R_+$ $\Phi'(B_{12}(x))$ commutes with $\Phi(B_{12}(1))$, with
$S_{(i,i+1)}$ for $i=3,\dots,n-1$, and with $diag [1,1,\mu_3,\dots,\mu_n]$ for $\mu_3,\dots,
\mu_n\in R_+^*$, we have
$$
\Phi'(B_{12}(x))=\begin{pmatrix}
a(x)& b(x)& &&\\
0& a(x)&&&\\
&&d(x)&&\\
&&&\ddots&\\
&&&&d(x)
\end{pmatrix},\text{ where } a(x), b(x)\in R_+, d(x)\in Z_+^*(R^*).
$$
Now we will use the condition
\begin{multline*}
(B_{12}(x))^2=diag [2,1,\dots,1] B_{12}(x)diag [1/2,1,\dots,1]\Rightarrow\\
\begin{pmatrix}
a(x)^2& a(x)b(x)+b(x)a(x)& &&\\
0& a(x)^2&&&\\
&&d(x)^2&&\\
&&&\ddots&\\
&&&&d(x)^2
\end{pmatrix}=
\begin{pmatrix}
a(x)& 2b(x)& &&\\
0& a(x)&&&\\
&&d(x)&&\\
&&&\ddots&\\
&&&&d(x)
\end{pmatrix}.
\end{multline*}
Consequently $a(x)=d(x)=1$ for every $x\in R_+$.
Therefore, if $\gamma=0$, then $\Phi'(B_{12}(x))=B_{12}(b(x))$ for every 
$x\in R_+$. Similarly, in the case $\beta=0$ we have $\Phi'(B_{12}(x))=B_{21}
(c(x))$ for every $x\in R_+$.
\end{proof}

Now the cases $\gamma=0$ and $\beta=0$ will be considered separately. 
In the next lemma we will prove that the case $\beta=0$ is impossible.

\begin{lemma}
If $n\ge 3$, $1/2\in R$, the automorphism
$\Phi'\in Aut(G_n(R))$ is such that $\forall 
\sigma\in \Sigma_n$ $\Phi'(S_\sigma)=S_\sigma$, then the condition
$$
\Phi'(B_{12}(1))=B_{21}(c(1))\quad (=B_{21}(\gamma))
$$
is impossible.
\end{lemma}
\begin{proof}
If $\Phi'(B_{12}(1))=B_{21}(\gamma)=B_{21}(c(1))$ for some $\gamma=c(1)\in R_+$, then
by the previous lemma there exists such a mapping 
$c(\cdot): R_+\to R_+$, that for every $x\in R_+$ we have $\Phi'(B_{12}(x))=B_{21}
(c(x))$. Since $n\ge 3$ and $\forall \sigma\in \Sigma_n$ $\Phi'(S_\sigma)=S_\sigma$,
we have 
$$
\Phi' (B_{13}(x))=\Phi'(S_{(2,3)} B_{12}(x)S_{(2,3)})=S_{(2,3)}B_{21}(c(x))S_{(2,3)}=
B_{31}(c(x)).
$$
Similarly, $\Phi' (B_{32}(x))=B_{23}(c(x))$.

Let us use the condition 
$$
\forall x_1,x_2\in R_+\quad
B_{13}(x_1)B_{32}(x_2)=B_{32}(x_2)B_{13}(x_1) B_{12}(x_1x_2).
$$
It implies
\begin{multline*}
\Phi'(B_{13}(x_1)B_{32}(x_2))=
\Phi'(B_{32}(x_2)B_{13}(x_1)B_{12}(x_1x_2))\Rightarrow\\
\Rightarrow B_{31}(c(x_1))B_{23}(c(x_2))=B_{23}(c(x_2))B_{31}(c(x_1))
B_{21}(c(x_1x_2))\Rightarrow\\
\Rightarrow
\begin{pmatrix}
1& 0& 0\\
0& 1& c(x_2)\\
c(x_1)& 0& 1
\end{pmatrix}=
\begin{pmatrix}
1& 0& 0\\
c(x_2)c(x_1)+c(x_1x_2)& 1& c(x_2)\\
c(x_1)& 0& 1
\end{pmatrix}\Rightarrow\\
\Rightarrow \forall x_1,x_2\in R_+\quad c(x_2)c(x_1)+ c(x_1x_2)=0\Rightarrow\\
\Rightarrow  \forall x\in R_+\quad c(x)^2+c(x^2)=0\Rightarrow c(x)=0,
\end{multline*}
but it is impossible, since $\Phi'$ is an autmorphism.
\end{proof}

Recall (see Definition~11) that if $G$ is some semigroup,
then a 
homomorphism $\lambda(\cdot ): G\to G$ is called a\emph{ central homomorphism}
of~$G$, 
if $\lambda(G)\subset Z(G)$. The mapping $\Omega(\cdot): G\to G$ such that 
$\forall X\in G$
$$
\Omega (X)=\lambda(X)\cdot X,
$$
where $\lambda(\cdot)$ is a central homomorphism, is called a \emph{central
homothety}.

Recall that for every $y(\cdot)\in Aut (R_+)$ by $\Phi^y$ we denote the automorphism 
of the semigroups $G_n(R)$ such that $\forall X=(x_{ij})\in G_n(R)$
$\Phi^y(X)=\Phi^y((x_{ij}))=(y(x_{ij}))$.

\begin{theorem}
Let $\Phi$ be any automorphism of $G_n(R)$, $n\ge 3$, $1/2\in R$. 
Then on the semigroup $GE_n^+(R)$ \emph{(}see Definition~\emph{10)} 
$\Phi=\Phi_M\Phi^c\Omega$, where 
$M\in \Gamma_n(R)$, $c(\cdot)\in Aut(R_+)$, $\Omega(\cdot)$ is a central homothety of
$GE_n^+(R)$.
\end{theorem}
\begin{proof}
By Lemma~6 there exists such a matrix $M'\in \Gamma_n(R)$ that for every $\sigma\in
\Sigma_n$ 
$$
\Phi'(S_\sigma)=\Phi_{M'}\Phi(S_\sigma)=S_\sigma.
$$

Consider now the automorphism~$\Phi'$. 

By Lemmas~12 and~13 there exists a mapping $c(\cdot): R_+\to R_+$ such that for every
$x\in R_+$
$$
\Phi'(B_{12}(x))=B_{12}(c(x)).
$$

Let us consider this mapping. Since $\Phi$ is an automorphism of the semigroup
$G_n(R)$, we have that $c(\cdot):R_+\to R_+$ is bijective.

Since for all $x_1,x_2\in R_+$\  $B_{12}(x_1+x_2)=B_{12}(x_1)B_{12}(x_2)$, we have
$B_{12}(c(x_1+x_2))=\Phi'(B_{12}(x_1+x_2))=\Phi'(B_{12}(x_1)
B_{12}(x_2))=\Phi'(B_{12}(x_1))\Phi'(B_{12}(x_2))=B_{12}(c(x_1))\cdot 
B_{12}(c(x_2))=B_{12}(c(x_1)+c(x_2))$, therefore for all $x_1,x_2\in R_+$ $c(
x_1+x_2)=c(x_1)+c(x_2)$ and $c(\cdot)$ is additive.

To prove the multiplicativity  of $c(\cdot)$ we will use the following:

1) $\Phi'(B_{13}(x))=\Phi'(S_{(2,3)}B_{12}(x)S_{(2,3)})=S_{(2,3)})=S_{(2,3)}B_{12}
(c(x))S_{(2,3)}=B_{13}(c(x))$;

2) similarly, $\Phi'(B_{32}(x))=B_{32}(c(x))$;

3) (compare with the proof of Lemma~13)
\begin{multline*}
B_{13}(x_1)B_{32}(x_2)=B_{32}(x_2)B_{13}(x_1)B_{12}(x_1x_2)\Rightarrow\\
\Rightarrow \Phi'(B_{13}(x_1))\Phi'(B_{32}(x_2))=\Phi'(B_{32}(x_2))\Phi'
(B_{13}(x_1))\Phi'(B_{12}(x_1x_2))\Rightarrow\\
\Rightarrow B_{13}(c(x_1))B_{32}(c(x_2))=B_{32}(c(x_2))B_{13}(c(x_1))B_{12}
(c(x_1x_2))\Rightarrow\\
\Rightarrow \forall x_1,x_2\in R_+\
\begin{pmatrix} 1& c(x_1)c(x_2)& c(x_1)\\
0& 1& 0\\
0& c(x_2)& 1
\end{pmatrix}=\begin{pmatrix}
1& c(x_1x_2)& c(x_1)\\
0& 1& 0\\
0& c(x_2)& 1
\end{pmatrix}\Rightarrow\\
\Rightarrow \forall x_1,x_2\in R_+\ c(x_1x_2)=c(x_1)c(x_2).
\end{multline*}

Therefore $c(\cdot)$ is a multiplicative mapping.

Since $c(\cdot)$ is a bijective, additive and multiplicative mapping, we have that
$c(\cdot )$ is an automorphism of the semiring~$R_+$, or, in other words, $c(\cdot)$ 
can be extended to an automorphism of the ring~$R$, preserving the order.

Consider now the mapping $\Phi^{c^{-1}}$, which maps every matrix $A=(a_{ij})$ to the matrix 
$\Phi^{c^{-1}}(A)=(c^{-1}(a_{ij}))$. This mapping is an automorphism of the semigroup
$G_n(R)$. Then $\Phi''=\Phi^{c^{-1}}\circ \Phi'=\Phi^{c^{-1}}\circ \Phi_{M'} \circ \Phi$ 
is an automorphism of the semigroup $G_n(R)$, preserving all matrices $S_\sigma$ ($
\sigma \in \Sigma_n$) and $B_{ij}(x)$ ($x\in R_+$, $i,j=1,\dots, n$, $i\ne j$).
Namely, $\Phi''(S_\sigma)=\Phi^{c^{-1}}(\Phi'(S_\sigma))=\Phi^{c^{-1}}(S_\sigma)=
S_\sigma$, since the matrix $S_\sigma$ contains only $0$ and~$1$;
for $i=3,\dots,n$ $\Phi'' (B_{i2}(x))=\Phi''(S_{(1,i)}B_{12}(x)S_{(1,i)})=
S_{(1,i)}\Phi''(B_{12}(x)))S_{(1,i)}=S_{(1,i)}\Phi^{c^{-1}}(B_{12}(c(x)))
S_{(1,i)}=S_{(1,i)}B_{12}(x)S_{(1,i)}=B_{i,2}(x)$;
for $j=3,\dots ,n$ 
$\Phi''(B_{1j}(x))=\Phi''(S_{(2,j)} B_{12}(x)S_{(2,j)})=S_{(2,j)} B_{12}(x)
S_{(2,j)}=B_{1j}(x)$;
for $i,j=3,\dots,n$ 
$\Phi''(B_{ij}(x))=\Phi''(S_{(i,1)}B_{1j}(x)S_{(1,i)})=S_{(1,i)}B_{1j}(x)
S_{(1,i)}=B_{ij}(x)$.

As we know (see Lemma~8), for all $\alpha\in R_+^*$
$$
\Phi''(diag [\alpha,1,\dots,1])=diag [\beta(\alpha),\gamma(\alpha),\dots,\gamma
(\alpha)],\quad \beta,\gamma\in R_+^*.
$$
Let us use the condition 
\begin{multline*}
diag[\alpha,1,\dots,1]B_{12}(1)diag[\alpha^{-1},1,\dots,1]=B_{12}(\alpha)\Rightarrow\\
\Phi''(diag[\alpha,1,\dots,1])\Phi''(B_{12}(1))\Phi''(diag [\alpha^{-1}
,1,\dots,1])=\Phi''(B_{12}(\alpha))\Rightarrow\\
\Rightarrow diag[\beta(\alpha),\gamma(\alpha),\dots,\gamma(\alpha)]B_{12}(1)diag [\beta
(\alpha)^{-1},\gamma(\alpha)^{-1},\dots,\gamma(\alpha)^{-1}]=B_{12}(\alpha)\Rightarrow\\
\Rightarrow \beta(\alpha)\gamma(\alpha)^{-1}=\alpha\Rightarrow \beta(\alpha)=\alpha
\gamma(\alpha)\Rightarrow\\
\forall \alpha\in R_+^*\ \Phi''(diag[\alpha,1,\dots,1])=diag [\alpha \gamma(\alpha),
\gamma(\alpha),\dots,\gamma(\alpha)].
\end{multline*}
Since $diag [\alpha,1,\dots,1]$ commutes with every matrix
of the form 
$$
\begin{pmatrix}
1& 0\\
0& X
\end{pmatrix},\quad X\in G_{n-1}(R),
$$
and $n\ge 3$, we have that for all $\alpha\in R_+^*$ $\gamma(\alpha)\in Z_+^*(R)$.

Since for all $\alpha_1,\alpha_2\in R_+^*$
\begin{multline*}
diag[\alpha_1\alpha_2\gamma(\alpha_1\alpha_2),\gamma(\alpha_1\alpha_2),\dots,
\gamma(\alpha_1\alpha_2)]=\Phi''(diag [\alpha_1\alpha_2,1,\dots,1])=\\
=\Phi''(diag[\alpha_1,1,\dots,1])\Phi''(diag[\alpha_2,1,\dots,1])=\\
=diag[\alpha_1\gamma(\alpha_1),\gamma(\alpha_1),\dots,\gamma(\alpha_1)]diag
[\alpha_2\gamma(\alpha_2),\gamma(\alpha_2),\dots,\gamma(\alpha_2)]=\\
=diag [\alpha_1\alpha_2\gamma(\alpha_1)\gamma(\alpha_2),\gamma(\alpha_1)\gamma(\alpha_2),
\dots,\gamma(\alpha_1)\gamma(\alpha_2)]\Rightarrow\\
\Rightarrow \forall \alpha_1,\alpha_2\in R_+^*\ 
\gamma(\alpha_1\alpha_2)=\gamma(\alpha_1)\gamma(\alpha_2),
\end{multline*}
we have that the mapping $\gamma(\cdot)$ is a central homomorphism (see Definition~11)
$\gamma(\cdot): R_+^*\to Z_+^*(R)$. 

If $A=diag[\alpha_1,\dots,\alpha_n]\in D_n(R)$, then
\begin{multline*}
\Phi''(A)=\Phi''(diag[\alpha_1,1,\dots,1] S_{1,2}diag[\alpha_2,1,\dots,1] S_{(1,2)}
S_{(1,3)} diag [\alpha_3,1,\dots,1]S_{(1,3)} \dots\\
\dots S_{(1,n)} diag[\alpha_n,1,
\dots,1]S_{(1,n)})=\gamma(\alpha_1)diag[\alpha_1,1,\dots,1]S_{(1,2)}\gamma(\alpha_2) diag
[\alpha_2,1,\dots,1]S_{(1,2)}\dots\\
\dots S_{(1,n)} \gamma(\alpha_n) diag[\alpha_n,1,\dots,1]
\gamma(\alpha_n)=\\
=\gamma(\alpha_1)\dots \gamma(\alpha_n)A=\gamma(\alpha_1\dots\alpha_n)A.
\end{multline*}

Recall (see Definition~8) that $\mathbf P$ is the subsemigroup in $G_n(R)$, which is generated
by the matrices~$S_\sigma$ ($\sigma\in \Sigma_n$), $B_{ij}(x)$ ($x\in R_+$, $i,j=1,\dots,n$,
$i\ne j$), and $diag [\alpha_1,\dots,\alpha_n]$ ($\alpha_1,\dots,\alpha_n\in R^*_+$).

It is clear that every matrix $A\in \mathbf P$ can be represented in the form
$$
A=diag[\alpha_1,\dots,\alpha_n]A_1\dots A_k,
$$
where $\alpha_1,\dots,\alpha_n\in R^*_+$, 
$A_1,\dots,A_k\in \{ S_\sigma, B_{ij}(x)|
\sigma\in \Sigma_n, x\in R_+, i,j=1,\dots,n, i\ne j\}$. 
Then 
\begin{multline*}
\Phi''(A)=\Phi''(diag[\alpha_1,\dots,\alpha_n]A_1\dots A_k)=\\
=\gamma(\alpha_1 \dots \alpha_n)diag [\alpha_1,\dots,\alpha_n]A_1\dots A_k=
\gamma (\alpha_1\dots \alpha_n)A.
\end{multline*}
Now we can introduce the mapping $\overline \gamma(\cdot): {\mathbf P}\to Z_+^*(R)$
by the following rule: if $A\in \mathbf P$ and $A=diag[\alpha_1,\dots,\alpha_n]A_1\dots A_k$, 
where 
$A_1,\dots,A_k\in \{ S_\sigma, B_{ij}(x)| \sigma\in \Sigma_n, x\in R_+, i,j=1,\dots,n,
i\ne j\}$, then $\overline \gamma (A)=\gamma(\alpha_1,\dots,\alpha_n)$.

The mapping $\overline \lambda(\cdot)$ is uniquely defined, because if 
$A=diag[\alpha_1,\dots,\alpha_n]A_1\dots A_k=diag [\alpha_1', \dots, \alpha_n']A_1'\dots
A_m'$, then $\Phi''(A)=\gamma(\alpha_1\dots \alpha_n)A$ and $\Phi''(A)=\gamma(\alpha_1'
\dots \alpha_n')A$ and therefore $\gamma(\alpha_1\dots \alpha_n)=\gamma(\alpha_1'\dots\alpha_n')$.

Since $\overline \gamma (AA') AA'=\Phi''(AA')=\Phi''(A)\Phi''(A')=\overline \gamma
(A) A\cdot \overline \gamma(A')A'=\overline \gamma(A) \overline \gamma(A')AA'$, 
we have that $\overline \gamma$ is a homomorphism $\mathbf P\to Z_+^*(R)$.

Now we see that on the semigroup~$\mathbf P$ the automorphism $\Phi''$ concides 
with the central homothety $\Omega(\cdot): \mathbf P\to \mathbf P$, where for all $a\in 
\mathbf P$ $\Omega(A)=\overline \gamma(A)\cdot A$. 

Let $B\in GE_n^+(R)$. Then (see Definitions~9,10) $B$ is $\cal P$-equivalent
to some matrix $A\in \mathbf P$, i.e. there exist matrices $A_0,\dots, A_k\in G_n(R)$,
$A_0=A\in \mathbf P$, $A_k=B$ and matrices $P_i,\widetilde P_i, Q_i, \widetilde Q_i\in 
\mathbf P$, $i=0,\dots,k-1$ such that for all $i=0,\dots,k-1$
$$
P_i A_i \widetilde P_i=Q_i A_{i+1}\widetilde Q_i.
$$
Then 
\begin{multline*}
\Phi''(P_0A_0\widetilde P_0)=\Phi''(Q_0A_1\widetilde Q_0)\Rightarrow\\
\Rightarrow \overline \gamma (P_0)P_0\overline \gamma (A_0)A_0\overline \gamma (\widetilde P_0)\widetilde 
P_0=\overline \gamma(Q_0)Q_0\Phi''(A_1)\overline \gamma(\widetilde Q_0)\widetilde Q_0
\Rightarrow\\
\overline \gamma (P_0A_0\widetilde P_0)P_0A_0\widetilde P_0=\overline \gamma (Q_0\widetilde 
Q_0)Q_0 \Phi''(A_1)\widetilde Q_0\Rightarrow\\
\Rightarrow \overline \gamma (P_0A_0\widetilde P_0)\overline \gamma(Q_0\widetilde 
Q_0)^{-1} Q_0A_1\widetilde Q_0=Q_0\Phi''(A_1)\widetilde Q_0\Rightarrow\\
\Rightarrow \Phi''(A_1) =\overline \gamma (P_0A_0\widetilde P_0)\overline \gamma
(Q_0\widetilde Q_0)^{-1}A_1,\dots,\\
\dots, \Phi''(B)=\Phi''(A_n)=\overline \gamma (P_{n-1})\overline \gamma (A_{n-1})
\overline \gamma (\widetilde P_{n-1})\overline \gamma (Q_{n-1})^{-1} \overline 
\gamma (\widetilde Q_{n-1}).
\end{multline*}
Therefore we can continue the mapping $\overline \gamma(\cdot): \mathbf P\to Z_+^*(R)$
to some mapping $\lambda(\cdot): GE_n^+(R)\to Z_+^*(R)$ such that for every
$B\in GE_n^+(R)$
$$
\Phi''(B)=\lambda (B)\cdot B.
$$

Since $\Phi''$ is an automorphism of the semigroup $GE_n^+(R)$, we have that
$\lambda(\cdot)$ is a central homomorphism $\lambda(\cdot): GE_n^+(R)\to Z_+^*(R)$ and so 
the automorphism $\Phi'':GE_n^+(R)\to GE_n^+(R)$ is a central homothety
$\Omega(\cdot): GE_n^+(R)\to GE_n^+(R)$, where $\forall X\in GE_n^+(R)$ $\Omega(X)=\lambda
(X)\cdot X$.

Since $\Phi''=\Omega$ on $GE_n^+(R)$ and $\Phi''=\Phi^{c^{-1}}\circ \Phi_{M'}\circ \Phi$
on $G_n(R)$, then $\Phi=\Phi_M\circ \Phi^c \circ \Omega$ on $GE_n^+(R)$, 
where $M={M'}^{-1}$.
\end{proof}

\end{document}